\documentclass[12pt]{article}

\usepackage{amsmath,amsfonts}
\usepackage{amssymb}
\usepackage{amsthm}

\usepackage{times,txfonts}

\newtheorem{thm}{Theorem}[section]

\newtheorem{lemma}[thm]{Lemma}
\newtheorem{dfn}[thm]{Definition}
\newtheorem*{claim}{\it Claim}
\newtheorem{remark}[thm]{\it Remark}
\newtheorem{example}[thm]{\it Example}
\numberwithin{equation}{section}

\def\pf{\noindent{\it Proof.} \ }  
\def\qed{\hfill $\square$}

\def\id{{\rm id}}
\def\dim{{\rm dim}~}

\title{
Tropical representation of Weyl groups  
associated 
with certain rational varieties}

\author{Teruhisa Tsuda and Tomoyuki Takenawa}
\date{}
\begin{document}

\maketitle
\begin{center}
\textit{ 
Dedicated to Professor Kazuo Okamoto
on his sixtieth birthday}
\end{center}

\begin{abstract}
Starting from certain rational varieties blown-up from $({\mathbb P}^1)^N$,
we construct a tropical,
i.e., subtraction-free birational,
representation of Weyl groups as a group of pseudo isomorphisms of the varieties.
We develop an  algebro-geometric framework of $\tau$-functions as 
defining functions of exceptional divisors on the varieties.
In the case where the corresponding root system is of affine type,
our construction yields a class of (higher order) $q$-difference Painlev\'e equations 
and its algebraic degree grows quadratically.
\end{abstract}

\renewcommand{\thefootnote}{\fnsymbol{footnote}}

\section{Introduction}
The aim of the present work is to develop 
the theory of birational representation of Weyl groups 
associated with algebraic varieties.

At the beginning of the twentieth century, 
it was discovered by Coble and Kantor, 
and later by Du Val,
that certain types of Cremona transformations act on 
the configuration space of point sets \cite{cob}.
Let $X_{m,n}$ be the configuration space of $n$ points in
general position in  the projective space ${\mathbb P}^{m-1}$.
Then, 
the Weyl group $W(T_{2,m,n-m})$ corresponding to the Dynkin diagram 
$T_{2,m,n-m}$
(see Figure~\ref{fig:dynkin})
acts birationally on $X_{m,n}$ 
and is generated by permutations of $n$ points and 
the standard Cremona transformation with respect to each $m$ points.
An algebro-geometric  and modern
interpretation  of this theory 
is due to
Dolgachev and Ortland   
\cite{do};
they showed that
the Cremona action of $W(T_{2,m,n-m})$ 
induces a {\it pseudo isomorphism},
i.e., an isomorphism except for subvarieties of codimension two or higher,
between
certain rational varieties blown-up from ${\mathbb P}^{m-1}$
at generic $n$ points,
which they call  
 {\it generalized Del Pezzo varieties}.
It is worth mentioning that
if $(m,n)=(3,9)$, the affine Weyl group of type $E_8^{(1)}$ appears and 
its lattice part gives rise to 
an important discrete dynamical system, i.e., 
the elliptic-difference Painlev\'e equation \cite{s};
see also \cite{org, kmnoy, kmnoy06}.

\begin{figure}[h]
\label{fig:dynkin}
\begin{center}
\begin{picture}(400,160)
\thicklines

%%%%%%%%%%%%%%%
%\put(80,170){$T_{2,m,n-m}$}
%%%%%%%%%%%%%%%%%%%%%%%%%%%
\put(40,87){$\underbrace{ \qquad \qquad \quad}$}
\put(66,70){$m$}

\put(101,87){$\underbrace{ \qquad \qquad \quad}$}
\put(118,70){$n-m$}

\put(92,98){$\left\{ \begin{array}{l}  \\ \end{array} \right.$}
\put(84,98){$2$}
%%%%%%%%%%%%%%%%%%%%%%%%%%%

\put(40,90){\circle{4}}  
\put(42,90){\line(1,0){8}}
\put(53,90){$\ldots$}
\put(78,90){\line(-1,0){8}}
\put(80,90){\circle{4}}    
\put(98,90){\line(-1,0){16}}
\put(100,90){\circle{4}}   
\put(102,90){\line(1,0){16}}
\put(120,90){\circle{4}}
\put(122,90){\line(1,0){8}}
\put(133,90){$\ldots$}
\put(158,90){\line(-1,0){8}}
\put(160,90){\circle{4}}   
\put(100,92){\line(0,1){16}}
\put(100,110){\circle{4}} 
  
%%%%%%%
%%%%%%%

%%%%%%%%%%%%%%%
%\put(295,170){$T^{\boldsymbol k}_{\boldsymbol \ell}$}
%%%%%%%%%%%%%%%%
%%
\put(290,120){$\left\{ \begin{array}{l} \\ \\ \\ \\ \end{array} \right.$}
\put(278,120){${k_n}$}
\put(290,55){$\left\{ \begin{array}{l} \\ \\ \\ \\ \end{array} \right.$}
\put(278,55){${\ell_n}$}

\put(248,120){$\left\{ \begin{array}{l} \\ \\ \\ \\ \end{array} \right.$}
\put(227,120){${k_{n-1}}$}
\put(248,55){$\left\{ \begin{array}{l} \\ \\ \\ \\ \end{array} \right.$}
\put(227,55){${\ell_{n-1}}$}

\put(332,120){$\left\{ \begin{array}{l} \\ \\ \\ \\ \end{array} \right.$}
\put(311,120){${k_{n+1}}$}
\put(332,55){$\left\{ \begin{array}{l} \\ \\ \\ \\ \end{array} \right.$}
\put(311,55){${\ell_{n+1}}$}

\put(300, 151){\circle{4}} 
\put(300,141){\line(0,1){8}}
\put(298,128){$\vdots$}
\put(300,116){\line(0,1){8}}
\put(300, 114){\circle{4}} 
\put(300,92){\line(0,1){20}}
\put(300,90){\circle{4}}  
\put(300,88){\line(0,-1){20}}
\put(300, 66){\circle{4}}
\put(300,64){\line(0,-1){8}}
\put(298,43){$\vdots$}
\put(300,39){\line(0,-1){8}}
\put(300,29){\circle{4}}   

\put(342, 151){\circle{4}} 
\put(342,141){\line(0,1){8}}
\put(340,128){$\vdots$}
\put(342,116){\line(0,1){8}}
\put(342, 114){\circle{4}}
\put(342,92){\line(0,1){20}}

\put(342,88){\line(0,-1){20}}
\put(342, 66){\circle{4}} 
\put(342,64){\line(0,-1){8}}
\put(340,43){$\vdots$}
\put(342,39){\line(0,-1){8}}
\put(342,29){\circle{4}} 

\put(258, 151){\circle{4}} 
\put(258,141){\line(0,1){8}}
\put(256,128){$\vdots$}
\put(258,116){\line(0,1){8}}
\put(258, 114){\circle{4}}
\put(258,92){\line(0,1){20}}

\put(258,88){\line(0,-1){20}}
\put(258, 66){\circle{4}} 
\put(258,64){\line(0,-1){8}}
\put(256,43){$\vdots$}
\put(258,39){\line(0,-1){8}}
\put(258,29){\circle{4}} 

\put(302,90){\line(1,0){38}}
\put(342, 90){\circle{4}}
\put(344,90){\line(1,0){20}}  
\put(368,90){$\ldots$}

\put(298,90){\line(-1,0){38}}
\put(258, 90){\circle{4}}
\put(256,90){\line(-1,0){20}}
\put(220,90){$\ldots$}

%%%
\put(298,5){$n$}
\put(332,5){$n+1$}
\put(244,5){$n-1$}
\put(220,5){$\cdots$}
\put(368,5){$\cdots$}
%%%
\end{picture}
\end{center}
\caption{Dynkin diagrams $T_{2,m,n-m}$ and $T^{\boldsymbol k}_{\boldsymbol \ell}$}
\end{figure}

In this paper,
starting from a certain rational variety
blown-up from 
$({\mathbb P}^1)^N$ 
along appropriate subvarieties
that are not only point sets,
we construct a birational representation of Weyl groups
corresponding to the Dynkin diagram   
$T^{\boldsymbol k}_{\boldsymbol \ell}$.
Here $T^{\boldsymbol k}_{\boldsymbol \ell}$ 
refers to the graph given in Figure~\ref{fig:dynkin},
specified by a pair of sequences ${\boldsymbol k}=(k_1, \ldots,k_N), {\boldsymbol \ell}=(\ell_1,\ldots, \ell_N) \in ({\mathbb Z}_{>0})^N$.
It is remarkable that
$T^{\boldsymbol k}_{\boldsymbol \ell}$ includes 
all of the simply-laced affine cases
$A_{n}^{(1)}$, $D_n^{(1)}$ and $E_n^{(1)}$,
which are relevant to a class of higher-dimensional discrete dynamical systems of Painlev\'e type.
This representation of Weyl groups is {\it tropical},
i.e., 
given in terms of subtraction-free 
birational mappings \cite{kir}
and,
interestingly enough,
possesses a geometric framework of $\tau$-functions.

In  the next section,  
we begin with blowing-up $({\mathbb P}^1)^N$ along certain subvarieties of codimension three.
Let $X$ be the rational variety thus obtained.
These $X$'s constitute a family.
Applying a cohomological technique (cf. \cite{do}),
we construct
the root and coroot lattices 
of type $T^{\boldsymbol k}_{\boldsymbol \ell}$
included in the N\'eron-Severi bilattice 
$N(X) \cong (H^2(X,{\mathbb Z}),H_2(X,{\mathbb Z}))$.
The associated Weyl group 
$W=W(T^{\boldsymbol k}_{\boldsymbol \ell})$
naturally acts on $N(X)$ as isometries (Lemma~\ref{lem:w}).
We see that this linear action of  
$W$ on $N(X)$
leads to 
a birational representation of $W$ 
on the family of varieties itself 
as a group of pseudo isomorphisms
(Theorem~\ref{thm:f-rep}). 
An element of $W$ naturally induces 
an appropriate permutation among the set of
exceptional divisors on $X$,
as similar to the classical topic:
27 lines (or exceptional curves)
on a cubic surface and a Weyl group of type $E_6$.
In \S~\ref{sect:tau},
for the purpose of describing 
the action of $W$
at the level of defining polynomials of exceptional divisors,
we introduce a geometric framework of $\tau$-functions
(Definition~\ref{def:tau} and Theorems~\ref{thm:tau-rep} and \ref{thm:tau-F}).
One important advantage
of our $\tau$-functions is that
we can trace
the resulting value for any element $w \in W$ 
by using only the defining polynomials of suitable divisors, 
although it is generally difficult 
to compute iterations of rational 
mappings.
In particular, our representation in an affine case
provides a discrete dynamical system 
 arising from the lattice part of the affine Weyl group.
Such a discrete dynamical system 
is equipped with a set of commuting discrete time evolutions and its algebraic degree grows 
in the quadratic order. 
And it is regarded 
as a (higher order)
$q$-difference Painlev\'e equation;
in \S~\ref{sect:aff}, the
$A_{n}^{(1)}$ and $D_n^{(1)}$
cases 
are demonstrated as typical examples.
In \S~\ref{sect:char},
we briefly indicate
(from  a soliton-theoretic  point of view)
some interesting relationships between
$\tau$-functions and the {\it character polynomials} 
appearing in representation theory of 
classical groups,
i.e., 
the Schur functions or the universal characters.

\section{Tropical Weyl group actions on rational varieties}
\label{sect:trop}

In this section, we consider a certain rational variety $X$
blown-up from 
$({\mathbb P}^1)^{N}$
and construct appropriate 
root and coroot lattices included in the N\'eron-Severi bilattice $N(X)$.
Moreover,
the action of the corresponding Weyl group 
is realized 
as pseudo isomorphisms 
of $X$
and is shown to possess a tropical 
(or subtraction-free birational) 
representation.

\subsection{Rational variety and root system}

Let 
 $f=(f_1,f_2,\ldots,f_N)$
 denote 
 the inhomogeneous coordinates of  $({\mathbb P}^1)^{N}$
 where $N \geq 3$.
 Fix 
 a pair of sequences
 ${\boldsymbol k}=(k_1, \ldots,k_N)$ and
 ${\boldsymbol \ell}=(\ell_1,\ldots, \ell_N)$
 of positive integers. 
 Consider the following subvarieties: 
 \begin{align*}
C_n^i
&=\left\{ f_{n-1}=0, f_n=-u_n^i, f_{n+1}=\infty \right\},
\quad
i=1,\ldots,k_n,
\\
C_n^{-j}
&=\left\{ f_{n-1}=\infty, f_n={-1}/{v_n^{-j}}, f_{n+1}=0 \right\},
\quad
j=1,\ldots, \ell_n,
\end{align*}
for $n=1,2,\ldots,N$,
where $u_n^i$ and $v_n^{-j}$ are nonzero parameters
(note: superscripts are indices, not exponents).
Hereafter we regard the suffix $n$ of the coordinate $f_n$ as an element of
${\mathbb Z}/N{\mathbb Z}$, namely,
$f_{n+N}=f_n$.
Let $\epsilon: X  \to ({\mathbb P}^1)^{N}$ be the blowing-up
along
$\{C_n^i, C_n^{-j}\}$.
Since $X$ is a rational variety, we have 
\[H^2(X,{\mathbb Z})\cong
NS(X)=\bigoplus_{n=1}^N 
\left(
{\mathbb Z} H_n \oplus
\bigoplus_{i=1}^{k_n} {\mathbb Z} E_n^i
\oplus 
\bigoplus_{j=1}^{\ell_n}  {\mathbb Z} E_n^{-j}
\right),
\]
where $NS(X)$ is the N\'eron-Severi group of $X$;
we denote by
$H_n$ the divisor class of hyperplanes $\{f_n= {\rm const.}\}$
and by $E_n^i$ the class of exceptional divisors $\epsilon^{-1}(C_n^i)$.
 The Poincar\'e duality guarantees 
$H_2(X,{\mathbb Z}) \cong (H^2(X,{\mathbb Z}))^*$.
 We can choose a basis $\{h_n, e_n^i\}$ of $H_2(X,{\mathbb Z})$,
 where $h_n$ corresponds to a line of degree 
 $\stackrel{ \stackrel{n}{\smile} }{ (0,\ldots,0,1,0,\ldots,0) }$
 and $e_n^i$ to a line restricted in a fibre ($\cong {\mathbb P}^2$) of the exceptional divisor $\epsilon^{-1}(C_n^i)$.
Thus  the intersection pairing $\langle  \,  ,   \, \rangle:H^2(X, {\mathbb Z}) \times H_2(X,{\mathbb Z}) \to {\mathbb Z}$ 
is defined by
$\langle  H_m  ,   h_n \rangle = \delta_{m,n}$, 
$\langle  E_m^i  ,   e_n^j \rangle = -\delta_{m,n}\delta_{i,j}$
and
$\langle {\rm otherwise}\rangle=0$.

Introduce the root lattice $Q$ and coroot lattice ${\check Q}$
as follows:
\[
Q= \bigoplus_{n=1}^N \bigoplus_{-\ell_n+1\leq i \leq k_n-1} {\mathbb Z}\alpha_n^i 
\subset
 H^2(X,{\mathbb Z})
\quad 
{\rm and}
\quad
{\check Q}=
\bigoplus_{n=1}^N\bigoplus_{-\ell_n+1\leq i \leq k_n-1} {\mathbb Z} {\check \alpha}_n^i
\subset H_2(X,{\mathbb Z}),
\] 
where 
\[
\begin{array}{ll}
 \alpha_n^0 =H_n-E_n^1-E_n^{-1},&
 {\check \alpha}_n^0=h_{n-1}+h_{n+1}-e_n^1-e_n^{-1}, 
 \\
 \alpha_n^i =E_n^i-E_n^{i+1}, &{\check \alpha}_n^i=e_n^i-e_n^{i+1} \quad (i=1,\ldots,k_n-1),
 \\
 \alpha_n^{-j} =E_n^{-j}-E_n^{-j-1}, &{\check \alpha}_n^{-j}
=e_n^{-j}-e_n^{-j-1} 
\quad(j=1,\ldots,\ell_n-1).
\end{array}
\]
For instance, we have $\langle\alpha_n^i, {\check \alpha}_n^i \rangle=-2$,
$\langle\alpha_n^0, {\check \alpha}_{n\pm 1}^0 \rangle=
\langle\alpha_{n\pm 1}^0, {\check \alpha}_n^0 \rangle=1$,
etc.
 The Dynkin diagram of the canonical root basis forms 
 $T^{\boldsymbol k}_{\boldsymbol \ell}$:
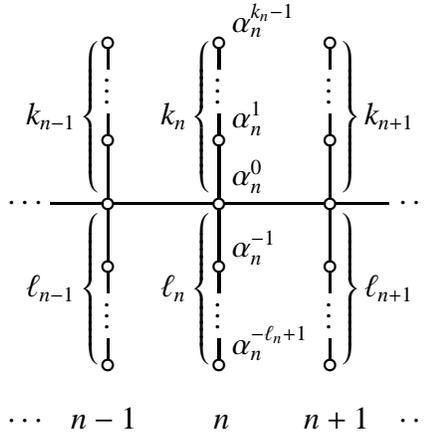
\begin{center}
\begin{picture}(200,170)
\thicklines
\put(90,120){$\left\{ \begin{array}{l} \\ \\ \\ \\ \end{array} \right.$}
\put(78,120){${k_n}$}

\put(90,55){$\left\{ \begin{array}{l} \\ \\ \\ \\ \end{array} \right.$}
\put(78,55){${\ell_n}$}

\put(48,120){$\left\{ \begin{array}{l} \\ \\ \\ \\ \end{array} \right.$}
\put(27,120){${k_{n-1}}$}

\put(48,55){$\left\{ \begin{array}{l} \\ \\ \\ \\ \end{array} \right.$}
\put(27,55){${\ell_{n-1}}$}

\put(135,120){$\left. \begin{array}{l} \\ \\ \\ \\ \end{array} \right\}$}
\put(155,120){${k_{n+1}}$}

\put(135,55){$\left. \begin{array}{l} \\ \\ \\ \\ \end{array} \right\}$}
\put(155,55){${\ell_{n+1}}$}

%%
%\put(101, 151){\circle{4}}  \put(105,155){$\alpha_{n}^{k_n-1}$}
\put(100, 151){\circle{4}}  \put(105,155){$\alpha_{n}^{k_n-1}$}
\put(100,141){\line(0,1){8}}
\put(98,128){$\vdots$}
\put(100,116){\line(0,1){8}}
\put(100, 114){\circle{4}} \put(105,119){$\alpha_{n}^1$}
\put(100,92){\line(0,1){20}}

\put(100,90){\circle{4}}  \put(105,97){$\alpha_{n}^0$}

\put(100,88){\line(0,-1){20}}
\put(100, 66){\circle{4}}  \put(105,70){$\alpha_{n}^{-1}$}
\put(100,64){\line(0,-1){8}}
\put(98,43){$\vdots$}
\put(100,39){\line(0,-1){8}}
\put(100,29){\circle{4}}   \put(105,33){$\alpha_{n}^{-\ell_n+1}$}

\put(142, 151){\circle{4}} 
\put(142,141){\line(0,1){8}}
\put(140,128){$\vdots$}
\put(142,116){\line(0,1){8}}
\put(142, 114){\circle{4}}
\put(142,92){\line(0,1){20}}

\put(142,88){\line(0,-1){20}}
\put(142, 66){\circle{4}} 
\put(142,64){\line(0,-1){8}}
\put(140,43){$\vdots$}
\put(142,39){\line(0,-1){8}}
\put(142,29){\circle{4}} 

\put(58, 151){\circle{4}} 
\put(58,141){\line(0,1){8}}
\put(56,128){$\vdots$}
\put(58,116){\line(0,1){8}}
\put(58, 114){\circle{4}}
\put(58,92){\line(0,1){20}}

\put(58,88){\line(0,-1){20}}
\put(58, 66){\circle{4}} 
\put(58,64){\line(0,-1){8}}
\put(56,43){$\vdots$}
\put(58,39){\line(0,-1){8}}
\put(58,29){\circle{4}}

\put(102,90){\line(1,0){38}}
\put(142, 90){\circle{4}}
\put(144,90){\line(1,0){20}}  
\put(168,90){$\ldots$}

\put(98,90){\line(-1,0){38}}
\put(58, 90){\circle{4}}
\put(56,90){\line(-1,0){20}}
\put(20,90){$\ldots$}

%%%
\put(98,5){$n$}
\put(132,5){$n+1$}
\put(44,5){$n-1$}

\put(20,5){$\cdots$}
\put(168,5){$\cdots$}

\end{picture}
\end{center}
The simple reflection $s_n^i$ 
associated with
a root $\alpha_n^i$ 
naturally acts on the N\'eron-Severi bilattice  
$N(X) \cong (H^2(X,{\mathbb Z}),H_2(X,{\mathbb Z}))$ 
as 
\begin{subequations}  \label{subeq:reflect}
\begin{align}
s_n^i(\Lambda)&=\Lambda + \langle  \Lambda , {\check \alpha}_n^i \rangle  \alpha_n^i,
\quad \Lambda \in H^2(X,{\mathbb Z}),
\\
s_n^i(\lambda)&=\lambda + \langle  \alpha_n^i, \lambda \rangle  {\check \alpha}_n^i,
\quad \lambda \in H_2(X,{\mathbb Z}).
\end{align}
\end{subequations}
One can easily check that these reflections indeed satisfy the 
fundamental relations (see \cite{kac})
of the Weyl group 
$W=W(T^{\boldsymbol k}_{\boldsymbol \ell})=\langle s_n^i \rangle$;
for instance, 
we have  $(s_n^i)^2=\id$ and $s_n^0 s_{n \pm 1}^0 s_n^0  = s_{n \pm 1}^0 s_n^0 s_{n \pm 1}^0$.

 The half of the anti-canonical class 
 $-\frac{1}{2}K_X= 
 \sum_{n=1}^N \left(H_n - \sum_{i=1}^{k_n} E_n^i-\sum_{j=1}^{\ell_n}E_n^{-j} \right) 
 \in H^2(X,{\mathbb Z}) $
 can be decomposed 
 in two ways 
 \[-\frac{1}{2}K_X=\sum_{n=1}^N D_n^0= \sum_{n=1}^ND_n^\infty,
 \]
 where
$D_n^0=H_n-\sum_{i=1}^{k_{n+1}} E_{n+1}^i-\sum_{j=1}^{\ell_{n-1}}E_{n-1}^{-j}$
and
$D_n^\infty=H_n-\sum_{i=1}^{k_{n-1}} E_{n-1}^i-\sum_{j=1}^{\ell_{n+1}}E_{n+1}^{-j}$.
Note that divisor classes $D_n^0$ and $D_n^\infty$ 
are effective and are represented by
the strict transforms of hyperplanes 
$\{f_n=0\}$ and $\{f_n=\infty\}$,
respectively.
In parallel, 
we shall {\it formally} define
an element 
$-\frac{1}{2}k_X\in H_2(X,{\mathbb Z})$
by
\[-\frac{1}{2}k_X=\sum_{n=1}^N \left(2h_n - \sum_{i=1}^{k_n} e_n^i-\sum_{j=1}^{\ell_n}e_n^{-j} \right) 
=\sum_{n=1}^N d_n^0= \sum_{n=1}^N d_n^\infty,
\]
where
$d_n^0=h_{n-1}+h_{n+1}-\sum_{i=1}^{k_{n+1}} e_{n+1}^i-\sum_{j=1}^{\ell_{n-1}}e_{n-1}^{-j}$
and
$d_n^\infty=h_{n-1}+h_{n+1}-\sum_{i=1}^{k_{n-1}} e_{n-1}^i-\sum_{j=1}^{\ell_{n+1}}e_{n+1}^{-j}$.
We see that $Q \subset \{d_n^0, d_n^\infty \}^\perp $ and
${\check Q}  \subset \{D_n^0, D_n^\infty\}^\perp$.
Hence we have the

\begin{lemma}  
\label{lem:w}
All elements $w \in W \subset {\rm Aut}(N(X))$ leave the intersection pairing
$\langle \, , \, \rangle$, $\frac{1}{2}K_X$ and $\frac{1}{2}k_X$ invariant.
\end{lemma}

\subsection{Birational representation of Weyl groups}

We summarize below the linear action of generators $s_n^i$ on the basis of $H^2(X, {\mathbb Z})$:
\begin{align*}
s_n^0(H_{n\pm1})&=H_{n\pm1}+H_n-E_n^1-E_n^{-1},
\quad
s_n^0(E_n^{\pm 1})=H_n-E_n^{\mp 1},
\\
s_n^i(E_n^{\{i,i+1\}})&=E_n^{\{i+1,i\}}  \quad {\rm for} \quad
1 \leq i \leq k_n-1,
\\
s_n^{-j}(E_n^{\{-j,-j-1\}})&=E_n^{\{-j-1,-j\}} \quad {\rm for} \quad 
1 \leq j \leq \ell_n-1.
\end{align*}
Now, let us extend the above linear action of the Weyl group 
$W=\langle s_n^i \rangle$ 
to the level of birational transformations 
on the rational variety $X$. 

To this end,
we first introduce the {\it multiplicative root variables} 
$a_n^i \in {\mathbb C}^\times$
attached to the canonical roots $\alpha_n^i$ and
fix the action of $s_n^i$ on them as 
\begin{equation}  \label{eq:mroot}
s_n^i(a_n^i)=\frac{1}{a_n^i},  \quad
s_n^i(a_n^{i  \pm 1})= a_n^i a_n^{i  \pm 1},  \quad
s_n^0(a_{n\pm1}^0)=a_n^0 a_{n\pm1}^0.
\end{equation}
Using the root variables, we fix the parameterization    
of subvarieties $C_n^i$ as follows:
\begin{equation}
u_n^1=u_n=
\frac{ \left( a_n^0  \prod_{j=1}^{\ell_n-1} (a_n^{-j})^{1- {j}/{\ell_n}} \right)^{{\ell_n}/{(k_n+\ell_n)}} }{
\left(\prod_{i=1}^{k_n-1} (a_n^{i})^{1- {i}/{k_n}} 
\right)^{{k_n}/{(k_n+\ell_n)} } },
\quad
v_n^{-1}=v_n
=
\frac{
\left(a_n^0\prod_{i=1}^{k_n-1} (a_n^{i})^{1- {i}/{k_n}} 
\right)^{{k_n}/{(k_n+\ell_n)} } }{ \left(\prod_{j=1}^{\ell_n-1} (a_n^{-j})^{1- {j}/{\ell_n}} \right)^{{\ell_n}/{(k_n+\ell_n)}} },
\end{equation}
and $u_n^i=a_n^{i-1}u_n^{i-1}$, 
$v_n^{-j}=a_n^{-j+1}v_n^{-j+1}$
for $i, j \geq 2$.
Namely, 
a projective equivalence class of the arrangement of subvarieties 
$\{C_n^i\}$ in $({\mathbb P}^1)^N$ can be identified with a point of 
the Cartan subalgebra 
associated with the Dynkin diagram $T^{\boldsymbol k}_{\boldsymbol \ell}$.
The rational varieties $X$'s under consideration constitute a family parameterized by 
the multiplicative root variables $a=(a_n^i)$;
so we shall write clearly as $X=X_{a}$.

Next, for each $w \in W \subset {\rm Aut}(H^2(X,{\mathbb Z}))$,
there is  a birational mapping ${\rm cr}(w): X_a \to X_{w(a)}$
such that $w={\rm cr}(w)^*$ (the pullback).
Hereafter we will intentionally omit to write clearly as ${\rm cr}(w)$,
and write shortly as $w$ instead.
Let 
${\mathbb K}(f)$ be
the field of rational functions in $f=(f_1,f_2,\ldots,f_N)$,
where the coefficient field ${\mathbb K}={\mathbb C}(a^{1/p})$ 
is generated by a suitable fractional power $(a_n^i)^{1/p}$ of $a_n^i$. 
We then have the following birational representation of 
$W$ over ${\mathbb K}(f)$,
acting on the family  ${\cal X}=\cup_a X_a$
 of varieties
as pseudo isomorphisms.

\begin{thm}  \label{thm:f-rep}
The birational transformations $s_n^i$ 
defined by
\begin{subequations}  \label{subeq:f-var}
\begin{align}
s_n^0(f_{n-1})
&=
(a_n^0)^{\ell_{n-1}/(k_{n-1}+\ell_{n-1})}
f_{n-1} 
\frac{f_n+  1/{v_n}}{f_n+ u_n},
\\
s_n^0(f_{n+1})
&=
(a_n^0)^{-k_{n+1}/(k_{n+1}+\ell_{n+1})}
f_{n+1}
\frac{f_n+  u_n }{f_n+  1/{v_n}},
\end{align}
and $s_n^i(f_m)=f_m$ (otherwise),
together with {\rm(\ref{eq:mroot})},
realize the Weyl group $W=W(T^{\boldsymbol k}_{\boldsymbol \ell})$ over the field ${\mathbb K}(f)$.
Moreover, each
$s_n^i$ maps $X_{a}$ to $X_{s_n^i(a)}$
as
a pseudo isomorphism of rational varieties 
and actually induces the linear action {\rm(\ref{subeq:reflect})} on $N(X)$.
\end{subequations}
\end{thm}

 \pf
 First we notice that
$s_n^0: u_n \leftrightarrow 1/v_n$,
$s_n^i: u_n^i \leftrightarrow u_n^{i+1}$ $(1 \leq i \leq k_n-1)$
and $s_n^{-j}: v_n^{-j}\leftrightarrow v_n^{-j-1}$ $(1 \leq j \leq \ell_n-1)$.
It is easy to verify by direct computation that
 $\langle s_n^i \rangle$ satisfies  the fundamental relations  of $W(T^{\boldsymbol k}_{\boldsymbol \ell})$.

For $i \neq 0$,
$s_n^i$ is trivial and  acts as  an isomorphism 
$X_a \to X_{s_n^i(a)}$,
so we have only to check that $s_n^0$ actually induces a pseudo isomorphism.
 The points of indeterminacy of the rational map 
 $s_n^0: ({\mathbb P}^1)^{N} \to  ({\mathbb P}^1)^{N}$
 are given by the following subvarieties:
 \[
 \begin{array}{ll}
 I_1(a)= \{f_n=-u_n, f_{n-1}=0 \},
 &  I_2(a)= \{f_n=-u_n, f_{n+1}=\infty \}, 
 \\
 I_3(a)= \{f_n=-1/{v_n}, f_{n+1}=0 \},
 &
  I_4(a)= \{f_n=-1/{v_n}, f_{n-1}=\infty \},
  \end{array}
 \]
 of codimension two.
 Since the same argument is also valid for the other $I_n$ $(n=2,3,4)$,
in what follows we shall only treat $I_1$.
We see that 
$s_n^0$ maps $I_1$ to $I_3$ generically, 
except for
$C_n^1=I_1 \cap I_2$;
furthermore,
$C_n^1$ equals
the inverse image
$(s_n^0)^{-1}(\{ f_n=-1/v_n\})$. 
By the blowing-up $\epsilon:X_a \to ({\mathbb P}^1)^{N}$,
the indeterminacy at $C_n^1$ is resolved and 
$s_n^0$ is extended holomorphically to $\epsilon^{-1}(C_n^1)$.
Thus we have proved that $s_n^0$ acts on $X_a$ as a pseudo isomorphism,
i.e., an isomorphism except for subvarieties of codimension two or higher.
 \qed

 \begin{remark} \rm
 Our representation given in Theorem~\ref{thm:f-rep} above
 is obviously tropical and hence admits a combinatorial counterpart 
via the ultra-discretization \cite{ttms}:
\[a \times b \to a +b, \quad a/b \to a-b,  \quad  a+b \to \min(a,b). \] 
 \end{remark}

\section{$\tau$-Functions}
\label{sect:tau}

An element of the Weyl group $W=W(T^{\boldsymbol k}_{\boldsymbol \ell})$
naturally induces a  permutation among the exceptional divisors 
on the rational variety $X$.
In this section, we introduce a geometric framework of 
{\it $\tau$-functions},
which describes the above permutations 
and therefore governs our tropical representation of 
$W$.

\subsection{Representation over the field of $\tau$-variables}

We set
\[
\theta_n^0= k_{n+1} +\ell_{n-1}
\quad {\rm and}  \quad
\theta_n^\infty=k_{n-1} +\ell_{n+1},
\]
which are equal to the numbers of centers of the blowing-up 
included in hyperplanes 
$\{f_n=0\}$ and  $\{f_n=\infty\}$, 
respectively.
From now on, 
we assume for simplicity that
\[
k_{n-1} k_{n+1}=\ell_{n-1} \ell_{n+1}.
\] 
In this case, 
(\ref{subeq:f-var})
can be equivalently rewritten into
\begin{subequations}  \label{subeq:f-var-2}
\begin{align}
s_n^0(f_{n-1})
&=
f_{n-1} 
\frac{ {v_n}^{\omega_n}f_n+ {v_n}^{-1+\omega_n} }{ {u_n}^{-\omega_n} f_n+ {u_n}^{1-\omega_n}},
\\
s_n^0(f_{n+1})
&=
f_{n+1}
\frac{ {u_n}^{-\omega_n} f_n+ {u_n}^{1-\omega_n}}{ {v_n}^{\omega_n}f_n+ {v_n}^{-1+\omega_n} },
\end{align}
\end{subequations} 
where a rational number $\omega_n$,
which is not a suffix,
is given by 
$\omega_n=\theta_n^0/(\theta_n^0+\theta_n^\infty)$.
Now we consider the decomposition of  variables 
\begin{equation}  \label{eq:f-x}
f_n=\frac{x_{n+1}}{x_{n-1}},
\end{equation}
and take a birational transformation
\begin{equation}  \label{eq:bir-x}
s_n^0(x_n)=x_n 
\frac{ {v_n}^{\omega_n}  x_{n+1} +{v_n}^{-1+\omega_n} x_{n-1}}{ {u_n}^{-\omega_n} x_{n+1} +{u_n}^{1-\omega_n} x_{n-1}}.
\end{equation}
One can easily verify that
(\ref{eq:bir-x}) with {\rm(\ref{eq:mroot})}
also defines a birational representation of   
$W$ on the field ${\mathbb K}(x)$.
Moreover,
the representation (\ref{subeq:f-var-2})
at the level of the $f$-variables 
is of course derived from (\ref{eq:bir-x})
via (\ref{eq:f-x}).
Introducing new variables $\tau_n^i$ 
attached to the centers $C_n^i$
of the blowing-up,
we
shall consider 
a further decomposition of variables in terms of the
$\tau$-variables
\begin{equation} \label{eq:x-tau} 
x_n=\frac{\xi_n}{\eta_n},   \quad
\xi_n= \prod_{i=1}^{k_n}\tau_n^i,
\quad 
\eta_n=\prod_{j=1}^{\ell_n}\tau_n^{-j}.
\end{equation}

\begin{thm}   \label{thm:tau-rep} 
Define the birational transformations $s_n^i$ by
\begin{subequations}    \label{subeq:tau-rep}
\begin{align}
\label{eq:tau-rep-a}
s_n^0(\tau_n^1)
&=
\frac{  {v_n}^{\omega_n}  \xi_{n+1}\eta_{n-1}+{v_n}^{-1+\omega_n} \xi_{n-1}\eta_{n+1}}{\tau_n^{-1} },
\\
\label{eq:tau-rep-b}
s_n^0(\tau_n^{-1})
&=
\frac{  {u_n}^{-\omega_n}  \xi_{n+1}\eta_{n-1}+{u_n}^{1-\omega_n} \xi_{n-1}\eta_{n+1}}{\tau_n^1 },
\\
s_n^i(\tau_n^{\{i,i+1\}})
&=\tau_n^{\{i+1,i\}} \quad (1 \leq i  \leq k_n-1),
\\
s_n^{-j}(\tau_n^{\{-j,-j-1\}})
&=\tau_n^{\{-j-1,-j\}}  \quad (1 \leq j  \leq \ell_n-1),
\\
s_n^i(\tau_m^j)
&=\tau_m^j   \quad  {\rm (otherwise)},
\end{align}
where 
$\omega_n=\theta_n^0/(\theta_n^0+\theta_n^\infty)$,
$\xi_n= \prod_{i=1}^{k_n}\tau_n^i$ and
$\eta_n=\prod_{j=1}^{\ell_n}\tau_n^{-j}$.
Then {\rm (\ref{subeq:tau-rep})} with {\rm(\ref{eq:mroot})} 
realizes the Weyl group
$W(T^{\boldsymbol k}_{\boldsymbol \ell})$ over the field 
$L={\mathbb K}(\tau)$
of rational functions
in indeterminates
$\tau_n^i$. 
\end{subequations}
\end{thm}

We remark that, 
conversely,
both realizations in terms of $f$-variables and $x$-variables 
are deduced  from Theorem~\ref{thm:tau-rep}
through
\begin{equation} \label{eq:f-tau}
f_n=\frac{\xi_{n+1} \eta_{n-1}}{\xi_{n-1} \eta_{n+1}},  
\quad
\xi_n= \prod_{i=1}^{k_n}\tau_n^i,
\quad 
\eta_n=\prod_{j=1}^{\ell_n}\tau_n^{-j}
\end{equation}
and through (\ref{eq:x-tau}), respectively.

\subsection{$\tau$-Function and its geometric meaning}

By means of Theorem~\ref{thm:tau-rep},
we now introduce the concept of 
{\it $\tau$-functions}; cf. \cite{kmnoy, t06}.
Let us consider the Weyl group orbit $M=W. \{ E_n^i\} \subset H^2(X,{\mathbb Z})$
of the classes of exceptional divisors.
Note that $\dim |\Lambda| =0$ for any $\Lambda \in M$.

\begin{dfn}  \label{def:tau}
\rm
The $\tau$-function
$\tau: M \to L$
is defined by
\\
\begin{tabular}{l}
(i) $\tau(E_n^i)=\tau_n^i$; \\ 
(ii) $\tau(w. \Lambda)= w. \tau(\Lambda)$ 
for any $\Lambda \in M$ and $w \in W$.
\end{tabular}
\end{dfn}

Here we suppose that an element $w \in W$
acts on a rational function $R(a,\tau) \in L$
as 
$w. R(a, \tau)=R(a. w, \tau. w)$,
that is, $w$ acts on the independent variables from the right.
By applying the Weyl group action given in Theorem~\ref{thm:tau-rep},
we can determine inductively
the explicit formula of 
$\tau(\Lambda)$ 
as a rational function in 
$\tau_n^i$.

Alternatively below
we shall present a geometric (and {\it a priori}) characterization of our $\tau$-functions.
Let us first
prepare some notations;
let $\zeta=(\zeta_1^0:\zeta_1^\infty,\ldots, \zeta_N^0 :\zeta_N^\infty)$
be the homogeneous coordinates of $({\mathbb P}^1)^{N}$
such that 
\begin{equation}
f_n=\frac{\zeta_n^0}{\zeta_n^\infty}.
\end{equation}
For a multiple index
$m=(m_1^0,m_1^\infty,\ldots, m_N^0,m_N^\infty) \in ({\mathbb Z}_{\geq 0})^{2N}$,
we set 
$\zeta^m=\prod_n (\zeta_n^0)^{m_n^0} (\zeta_n^\infty)^{m_n^\infty}$.
An element 
$\Lambda=\sum d_n H_n -\sum \mu_n^i E_n^i  \in M$
corresponds to a unique hypersurface of degree $d=(d_1,\ldots,d_N)$
passing through  
$C_n^i$ with multiplicity $\mu_n^i$.
Let 
\begin{equation}
\Phi_\Lambda (\zeta)=\sum_m  A_m \zeta^m
\in {\mathbb K}[\zeta]
\end{equation}
be the defining polynomial of the hypersurface,
where $A_m \in {\mathbb K}$ is a monic monomial in multiplicative root variables.
Normalize 
$\Phi_\Lambda$  
as its coefficients satisfy the following condition:
\begin{equation} \label{eq:norm}
\prod_m {A_m}^{ (1/\theta)^m}=1
\end{equation}
with
$(1/\theta)^m =\prod_n(1/\theta_n^0)^{m_n^0}(1/\theta_n^\infty)^{m_n^\infty}$.
A few examples of the normalized defining polynomial $\Phi_\Lambda(\zeta)$ $(\Lambda \in M)$
are 
\begin{align*}
\Phi_{E_n^i}
&=1, 
\\
\Phi_{H_n-E_n^1}
&=  {u_n}^{-\omega_n}  \zeta_n^\infty(f_n+u_n),
\\
\Phi_{H_n+H_{n+1}-E_n^1-E_n^{-1}-E_{n+1}^{-1}}
&={u_n}^{ \omega_n(\omega_{n+1}-1)}
{v_n}^{\omega_n \omega_{n+1} }
{v_{n+1}}^{\omega_{n+1}}
\\
& \times
\zeta_n^\infty \zeta_{n+1}^\infty
\left( f_n f_{n+1}+u_n f_{n+1}+\frac{f_n}{v_{n+1}}+\frac{1}{v_n v_{n+1}}  \right).
\end{align*}
Now assume that
\begin{equation}  \label{eq:sup}
\zeta_n^0=\xi_{n+1}\eta_{n-1},
\quad
\zeta_n^\infty=\xi_{n-1}\eta_{n+1}.
\end{equation}
We then have an expression of the $\tau$-function $\tau(\Lambda)$
in terms of the normalized defining polynomial 
$\Phi_\Lambda(\zeta)$
of the corresponding divisor $\Lambda$.

\begin{thm} \label{thm:tau-F}
For any 
$\Lambda=\sum_n d_n H_n -\sum_{n,i} \mu_n^i E_n^i  \in M$,
we have the equality
\begin{equation}   \label{eq:tau-F}
\tau(\Lambda)\prod_{n,i} (\tau_n^i)^{\mu_n^i}=\Phi_\Lambda(\zeta).
\end{equation}
\end{thm}

Hence,
one can  in principle trace 
the resulting value with respect to any iteration of Weyl group actions
by means of the defining polynomials of appropriate divisors.  
We see from (\ref{eq:tau-F}) that
$\tau(\Lambda)$ 
is a Laurent polynomial in the indeterminates $\tau_n^i$,
though, from the definition of $\tau$-functions,
we can only state that
$\tau(\Lambda)$ is rational.
This {\it regularity} is an interesting feature 
of our $\tau$-functions and
should be relevant to the theory of infinite-dimensional integrable systems 
such as the KP hierarchy, UC hierarchy, etc.;
see \S~\ref{sect:char}.  
Note that the $\tau$-function can be considered as a counterpart of height functions
in the sense that each of the original inhomogeneous coordinates 
$f=(f_1,f_2, \ldots,f_N)$ of $({\mathbb P}^1)^N$ is recovered as a ratio of $\tau$-functions;
see (\ref{eq:f-tau}). 
\\

\noindent
{\it Proof of Theorem~\ref{thm:tau-F}.}
We will prove the theorem by induction.
If
$\Lambda=H_n-E_n^{\pm 1}=s_n^0(E_n^{\mp 1}) \in M$,
then  
(\ref{eq:tau-F})
follows immediately
from (\ref{eq:tau-rep-a}) and (\ref{eq:tau-rep-b}).
We assume that  (\ref{eq:tau-F}) is true for 
$\Lambda  \in M$.
Then it is enough to verify for 
$\Lambda'=w(\Lambda)$
such that $w \in \{s_n^i \}$ is a generator of $W=W(T^{\boldsymbol k}_{\boldsymbol \ell})$.
The action of $s_n^i$ $(i \neq 0)$ 
is just a permutation of $\tau$-variables,
so we need only to
 concentrate upon the nontrivial case 
$w=s_n^0$.

For 
$\Lambda=\sum_n d_n H_n -\sum_{n,i} \mu_n^i E_n^i  \in M$,
let us introduce  
a polynomial $\varphi_\Lambda(f) \in {\mathbb K}[f]$
in the inhomogeneous coordinates 
$f=(f_1,\ldots,f_N)$ of degree $d=(d_1,\ldots,d_N)$
defined by
$\varphi_\Lambda(f)=\Phi_\Lambda(\zeta)\prod_{n}(\zeta_n^\infty)^{-d_n}$.

\begin{claim}\rm
We have
\begin{equation}  
s_n^0(\varphi_\Lambda(f)) 
=
c
(f_n+u_n)^{-d_{n-1}+\mu_n^{-1}}
(f_n+1/{v_n})^{-d_{n+1}+\mu_n^1}\varphi_{s_n^0.\Lambda}(f),
\label{eq:claim}
\end{equation}
where the normalizing constant $c$ is given by  
$c=(
{u_n}^{d_{n-1}-\mu_n^{-1}}
{v_n}^{-d_{n+1}+\mu_n^{1}}
)^{\omega_n}$.
\end{claim}

To prove this, first we notice that
$s_n^0 . \Lambda
=\Lambda+\langle \Lambda, \check \alpha_n^0\rangle \alpha_n^0
=\Lambda+(d_{n-1}+d_{n+1}-\mu_n^1-\mu_n^{-1})\alpha_n^0$;
accordingly
the degree $d'=(d_1',\ldots,d_N')$ of 
$\varphi_{s_n^0.\Lambda}(f)$
reads 
$d_n'=d_{n-1}+d_n+d_{n+1}-\mu_n^1-\mu_n^{-1}$
and
$d_i'=d_i$ $(i \neq n)$.
The multiplicity of 
$C_n^{1}=\{f_{n-1}=0,f_n=-u_n,f_{n+1}=\infty\}$ 
on
the hypersurface
$\{\varphi_\Lambda(f)=0\} \subset ({\mathbb P}^1)^{N}$
is given by
${\rm ord}_{C_n^1}(\varphi_\Lambda)=\mu_n^1$.
Hence the defining polynomial $\varphi_\Lambda(f)$ 
is of the form
\[
\varphi_\Lambda(f)
=\sum_{i \geq  0}\psi_{\mu_n^1-d_{n+1}+i},
\]
where $\psi_i=\psi_i(f)$ is homogeneous of degree $i$
in 
$(f_{n-1}, f_n+u_n, 1/{f_{n+1}})$.
In the same way, 
from ${\rm ord}_{C_n^{-1}}(\varphi_\Lambda)=\mu_n^{-1}$,
we can also describe it as 
$\varphi_\Lambda
=\sum_{i \geq  0}{\widetilde \psi}_{\mu_n^{-1}-d_{n-1}+i}$ where
${\widetilde \psi}_i$
is homogeneous of degree $i$ in 
$(1/{f_{n-1}},f_n+1/{v_n},f_{n+1})$.
By applying $(\ref{subeq:f-var-2})$,
we see that
$s_n^0(\varphi_\Lambda)
\times
(f_n+u_n)^{d_{n-1}-\mu_n^{-1}}
(f_n+1/{v_n})^{d_{n+1}-\mu_n^1} $
is a polynomial of degree at most 
$d'$
and therefore equals 
$({\rm const.}) \times \varphi_{s_n^0.\Lambda}$.
Here we used implicitly the fact that $s_n^0$ is a pseudo isomorphism
(Theorem~\ref{thm:f-rep})
and transforms $\{ \varphi_\Lambda(f) =0 \}$
generically  
into $\{ \varphi_{s_n^0. \Lambda} (f)=0 \}$.    
Since (\ref{subeq:f-var-2}) keeps the normalization (\ref{eq:norm}) invariant, 
we can immediately 
calculate the constant
and thus arrive at (\ref{eq:claim}).
\\

Applying $s_n^0$ to both sides of (\ref{eq:tau-F})
and using the claim above,
we obtain
\begin{align*}
({\rm LHS})
& \stackrel{s_n^0}{\longrightarrow}
\tau(\Lambda')\prod_{i,j}(\tau_j^i)^{\mu_j^i}
\times 
\left(  \frac{s_n^0(\tau_n^1)}{\tau_n^1} \right)^{\mu_n^1} 
\left(  \frac{s_n^0(\tau_n^{-1})}{\tau_n^{-1}} \right)^{\mu_n^{-1}},
\\
({\rm RHS})
&=\varphi_\Lambda(f) \prod_j(\zeta_j^\infty)^{d_j}
\\
&\stackrel{s_n^0}{\longrightarrow}
s_n^0(\varphi_\Lambda(f))\prod_j(\zeta_j^\infty)^{d_j}
\times 
\left(  \frac{s_n^0(\tau_n^1)}{\tau_n^1} \right)^{d_{n+1}} 
\left(  \frac{s_n^0(\tau_n^{-1})}{\tau_n^{-1}} \right)^{d_{n-1}}
\\
&=\varphi_{\Lambda'}(f)\prod_j(\zeta_j^\infty)^{d_j'} \times
\left(\frac{s_n^0(\tau_n^1)}{\tau_n^1} \right)^{\mu_n^1}\left(\frac{s_n^0(\tau_n^{-1})}{\tau_n^{-1}}\right)^{\mu_n^{-1}}
(\tau_n^1 \tau_n^{-1})^{-d_{n-1}-d_{n+1}+\mu_n^{1}+\mu_n^{-1}}.
\end{align*}
Here we have used the formulae
$f_n+u_n={u_n}^{\omega_n}s_n^0(\tau_n^{-1})\tau_n^1/\zeta_n^\infty$
and
$f_n+1/{v_n}={v_n}^{-\omega_n}s_n^0(\tau_n^{1})\tau_n^{-1}/\zeta_n^\infty$.
Thus, we get
\[
\Phi_{\Lambda'}(\zeta)=
\varphi_{\Lambda'}(f)\prod_j(\zeta_j^\infty)^{d_j'}
=\tau(\Lambda')\prod_{i,j}(\tau_j^i)^{\mu_j^i}  \times
(\tau_n^1 \tau_n^{-1})^{d_{n-1}+d_{n+1}-\mu_n^{1}-\mu_n^{-1}},
\]
which is exactly
(\ref{eq:tau-F}) for 
$\Lambda'=s_n^0 . \Lambda$.
The proof of the theorem is complete.
\qed

\begin{remark} \rm

Let us extendedly apply (\ref{eq:tau-F}) to the effective divisor classes 
$D_n^0=H_n-\sum_{i=1}^{k_{n+1}} E_{n+1}^i-\sum_{j=1}^{\ell_{n-1}} E_{n-1}^{-j}$
and
$D_n^\infty=H_n-\sum_{i=1}^{k_{n-1}} E_{n-1}^i-\sum_{j=1}^{\ell_{n+1}} E_{n+1}^{-j}$,
which correspond to  the hyperplanes
$\{\zeta_n^0=0\}$ and $\{\zeta_n^\infty=0\}$,
respectively; thus, 
we obtain 
\[\tau(D_n^0) \xi_{n+1}\eta_{n-1}=\zeta_n^0
\quad {\rm and} 
\quad
\tau(D_n^\infty) \xi_{n-1}\eta_{n+1}=\zeta_n^\infty.
\]
Because $D_n^0$ and $D_n^\infty$
are invariant under 
the action of
Weyl group $W=W(T^{\boldsymbol k}_{\boldsymbol \ell})$,
we can define
$\tau(D_n^0)=\tau(D_n^\infty) \equiv 1$,
which affirms the assumption (\ref{eq:sup}).
\end{remark}

\section{Examples of affine case and $q$-Painlev\'e equations}
\label{sect:aff}

We assume that the Dynkin diagram $T^{\boldsymbol k}_{\boldsymbol \ell}$ is of affine type.
Then the lattice part of the affine Weyl group
$W(T^{\boldsymbol k}_{\boldsymbol \ell})$
provides an interesting class of 
discrete dynamical systems of Painlev\'e type.
This is based on the following facts:
\\
{\bf (i) (Existence   of commuting time evolutions)}
For each root $\alpha \in Q$, 
we have the translation 
$t_{\alpha}$ acting on the N\'eron-Severi bilattice 
$N(X)\cong(H^2(X,{\mathbb Z}), H_2(X,{\mathbb Z}))$
by the Kac formula \cite[\S6.5]{kac}:
\begin{align}
t_\alpha(\Lambda) &= 
\Lambda - \langle \Lambda, \check{\delta} \rangle \alpha 
+ \left( \langle \Lambda, \check{\alpha} \rangle  
   - \frac{1}{2} \langle \alpha , \check{\alpha}\rangle
  \langle \Lambda,\check{\delta}  \rangle \right)\delta,
  \label{eq:kac}
\\
t_\alpha(\lambda) &=
\lambda - \langle   \delta ,\lambda \rangle \check{\alpha} 
+ \left( \langle \alpha, \lambda\rangle  
   - \frac{1}{2} \langle \alpha , \check{\alpha}\rangle
  \langle \delta, \lambda  \rangle \right)\check{\delta}
\end{align}
for $\Lambda\in H^2(X,{\mathbb Z})$ and $\lambda\in H_2(X,{\mathbb Z})$.
Here 
$\delta=-\frac{1}{2}K_X$ and $\check{\delta}=-\frac{1}{2}k_X$
are
$W(T^{\boldsymbol k}_{\boldsymbol \ell})$-invariants (Lemma~\ref{lem:w})
and play the role of the null vector and its dual,
respectively.
The additivity property of translations
$ t_\alpha t_\beta =t_\beta t_\alpha =t_{\alpha+\beta}$ holds.
This implies that the discrete dynamical system
\[t_\alpha(f_i)= F_{\alpha,i}(a,f)   \quad (i=1, \ldots, N)\]
corresponding to a root $\alpha \in Q$
has a set of commuting discrete time evolutions 
$\{ t_\beta \,  | \,  \beta \in Q \}$,
where $F_{\alpha,i}(a,f) \in {\mathbb K}(f)$ 
is a rational function explicitly determined by 
Theorem~\ref{thm:f-rep}; cf. \cite{ny}.
\\
{\bf (ii) (Quadratic degree growth)}
From (i), we have  
\begin{align*}
(t_{\alpha})^n(H_i)&=t_{n \alpha}(H_i)
\\
&= H_i - n \langle H_i, \check{\delta} \rangle \alpha 
+ \left( n \langle  H_i, \check{\alpha} \rangle  
   - \frac{n^2}{2} \langle \alpha , \check{\alpha}\rangle
  \langle H_i,\check{\delta}  \rangle \right)\delta.
\end{align*}
Therefore the value of pairing 
$\langle  (t_{\alpha})^n(H_i), h_j  \rangle$
is in the order $n^2$, that is, the degree of $(t_\alpha)^n(f_i)$ with respect to $f_j$ 
is at most $n^2$.

Generally,
a quadratic degree growth asserts integrability of the dynamical system under consideration;
in  \cite{tak01},
it was proved by a similar technique as above that
the degree growth of
every {\it second order} difference Painlev\'e equation (in Sakai's classification \cite{s})
is quadratic.
Also, as shown in \cite{tak04}, 
a certain class of integrable higher dimensional dynamics defined on 
generalized Del Pezzo varieties has a quadratic degree growth. 
More elementarily,  the degree of
the $n$-th iteration of an addition formula of an elliptic curve grows 
quadratically
as a rational mapping;
see, e.g., \cite{sil86}.

In this section, 
we shall demonstrate  the
$A_r^{(1)}$ and $D_r^{(1)}$
cases
as examples.

\subsection{The case of $A_{N-1}^{(1)}$}

Let us consider the case where
${\boldsymbol k}={\boldsymbol \ell}=(1,1,\ldots,1) \in ({\mathbb Z}_{>0})^N$
so that the corresponding
 Dynkin diagram $T^{\boldsymbol k}_{\boldsymbol \ell}$ is of type $A_{N-1}^{(1)}$.
We have,  
from
Theorems \ref{thm:f-rep} and \ref{thm:tau-rep}, 
a tropical representation of affine Weyl group 
$W(A_{N-1}^{(1)})$
as pseudo isomorphisms of a rational variety $X_a$.
Moreover,
we can extend  $W(A_{N-1}^{(1)})$
to  $W(A_{N-1}^{(1)})\times W(A_{1}^{(1)})$
in the following manner.

First  we reformulate the parameters $u_n$ and $v_n$ as follows:
\begin{equation}
u_n=\frac{a_n b_1}{b_0},
\quad
v_n=\frac{a_n b_0}{b_1},
\end{equation}
where $a_n$ 
and $b_i$ $(i=0,1)$ are nonzero parameters such that 
\begin{equation}
\prod_{n=1}^N a_n= (b_0 b_1)^N=q^N .
\end{equation}
Here $a_n$ and $b_i$ will play roles of multiplicative root variables 
of  $A_{N-1}^{(1)}$ and $A_{1}^{(1)}$,
respectively.

Let $\gamma_{1,n}$ and $\check{\gamma}_{1,n}$
be a pair of vectors defined by
\begin{align*}
\gamma_{1,n}&= -\frac{K_X}{2}-H_n-E_{n}^{1}+E_{n-1}^{-1}+E_{n}^{-1}+E_{n+1}^{-1}
\in H^2(X,{\mathbb Z}),
\\
\check{\gamma}_{1,n}&=h_n-e_n^{1}
\in H_2(X,{\mathbb Z})
\end{align*}
for $n \in {\mathbb Z}/N{\mathbb Z}$.
We see that these are mutually 
orthogonal
$(-2)$-vectors,
namely, 
$\langle \gamma_{1,i},\check{\gamma}_{1,j}\rangle=-2 \delta_{i,j}$.
Define the action of the reflection $r_{1,n}$ associated with
the pair
$(\gamma_{1,n}, \check{\gamma}_{1,n})$ 
on $N(X) \cong (H^2(X,{\mathbb Z}),H_2(X,{\mathbb Z}))$
in the same way as (\ref{subeq:reflect}).
Let 
$\beta_1=\sum_{n=1}^N \gamma_{1,n}$
and 
$\check{\beta}_1=\sum_{n=1}^N \check{\gamma}_{1,n}$
and define the associated reflection $r_1$ by
$r_1=r_{1,1} \circ \cdots \circ r_{1,N}$.
In addition, 
we take an involution 
$\iota: E_n^{1} \leftrightarrow E_n^{-1},
e_n^{1} \leftrightarrow e_n^{-1}$
of $N(X)$;
let $(\beta_{0},\check{\beta}_{0})=\iota(\beta_1,\check{\beta}_1)$
and define $r_0=\iota \circ r_1 \circ \iota$.
Then we see that $\langle r_0,r_1\rangle\cong W(A_1^{(1)})$
and $\iota$ realizes
the Dynkin diagram automorphism of $A_1^{(1)}$.
Likewise, we can introduce a diagram automorphism $\pi$ of $A_{N-1}^{(1)}$
 acting on $N(X)$ as
$\pi: H_n \mapsto H_{n+1}, E_{n}^{\pm 1} \mapsto E_{n+1}^{\pm 1}, h_n \mapsto h_{n+1}, e_{n}^{\pm 1} \mapsto e_{n+1}^{\pm 1}$.
Note that $(\beta_i, \check{\beta}_i)$ $(i = 0,1)$
is orthogonal to 
$(\alpha_n, \check{\alpha}_n)$,
where 
$\alpha_n=\alpha_n^0=H_n-E_n^1-E_n^{-1}$
and 
$\check{\alpha}_n=\check{\alpha}_n^0=h_{n-1}+h_{n+1}-e_n^1-e_n^{-1}$
are the (co-)root bases of $A_{N-1}^{(1)}$.
Furthermore, we have in fact
\begin{align*}
Q \oplus {\mathbb Z} \beta_1
&= 
\left(\bigoplus_{n=1}^N {\mathbb Z} \alpha_n\right)
\oplus {\mathbb Z} \beta_1
=\{d_n^0, d_n^\infty \}^\perp  \subset H^2(X,{\mathbb Z}),
\\
\check{Q} \oplus {\mathbb Z} \check{\beta}_1
&= 
\left(\bigoplus_{n=1}^N {\mathbb Z} \check{\alpha}_n\right)
\oplus {\mathbb Z} \check{\beta}_1
=\{D_n^0, D_n^\infty \}^\perp  \subset H_2(X,{\mathbb Z}).
\end{align*}

In order to realize $r_i$ $(i=0,1)$ as birational transformations,
we give the action of $r_1$ on the cohomology basis of $H^2(X,{\mathbb Z})$:
\[
r_{1}(E_n^{1})=
E_n^1+\gamma_{1,n}=
\sum_{i \neq n} H_i -\sum_{j=1}^N E_j^1 -\sum_{k \neq n,n\pm1}E_k^{-1},
\]
$r_{1}(E_n^{-1})=E_n^{-1}$
and
$r_{1}(H_n)=
H_n+\gamma_{1,n}$.
For $\Lambda=r_{1}(E_n^{1})$,
the normalized defining polynomial reads
\[
\Phi_\Lambda(\zeta)
=G_n(\zeta)
:=\left(
1+ \sum_{j=1}^{N-1} \prod_{i=1}^{j}\frac{u_{n+i}}{f_{n+i}}
\right)
\prod_{k=1}^{N-1}
(u_{n+k})^{-1+k/N}
\zeta_{n+k}^{0}.
\]
By virtue of Theorem~\ref{thm:tau-F},
we thus obtain
a birational action of $r_1$ on $\tau$-variables: 
\[r_1(\tau_n^1)
=\frac{G_n(\zeta)}{\prod_{j}\tau_j^1 \prod_{k \neq n,n\pm 1} \tau_k^{-1}},
\quad
r_1(\tau_n^{-1})=\tau_n^{-1}.
\]
Recall here that 
$\zeta_n^0=\tau_{n+1}^1\tau_{n-1}^{-1}$,
$\zeta_n^\infty=\tau_{n-1}^1\tau_{n+1}^{-1}$
and 
$f_n=\zeta_n^0/\zeta_n^\infty$.

Now we fix the action on 
variables $a_n$ 
and $b_i$
as 
\begin{equation}  \label{eq:para-a}
\begin{array}{l}
s_n(a_n)=1/{a_n}, \quad s_n(a_{n\pm 1}) =a_na_{n\pm 1}, \\
r_i(b_i)=1/{b_i}, \quad  r_i(b_j)={b_i}^2b_j  \quad (i \neq j), \\
\pi(a_n)=a_{n+1}, \quad \iota(b_{\{0,1\}}) =b_{\{1,0\}}.
\end{array}
\end{equation}
The following birational transformations 
$\langle s_n, \pi ,  r_i, \iota \rangle$ 
together with (\ref{eq:para-a})
realize the extended affine Weyl group
$W=\widetilde{W}(A_{N-1}^{(1)}) \times \widetilde{W}(A_{1}^{(1)})$:
\begin{subequations}  \label{subeq:weyla}
\begin{align}
s_n(\tau_n^{1})
&= \frac{{v_n}^{1/2}\tau_{n+1}^{1}\tau_{n-1}^{-1} +{v_n}^{-1/2}\tau_{n-1}^{1}\tau_{n+1}^{-1} }{\tau_n^{-1}},
\\
s_n(\tau_n^{-1})
&= \frac{{u_n}^{-1/2}\tau_{n+1}^{1}\tau_{n-1}^{-1} +{u_n}^{1/2}\tau_{n-1}^{1}\tau_{n+1}^{-1} }{\tau_n^{1}},
\\
r_1(\tau_n^1)
&=\frac{G_n(\zeta)}{\prod_{j} \tau_j^1 \prod_{k \neq n,n\pm 1}\tau_k^{-1}},
\\
\pi(\tau_n^{\pm 1})
&=\tau_{n+1}^{\pm 1}, \quad
\iota(\tau_n^{\pm 1})=\tau_{n}^{\mp 1},
\end{align}
\end{subequations}
and
$r_{0}= \iota \circ r_1 \circ\iota$.
Consequently, with respect to variables $f_n=(\tau_{n+1}^1\tau_{n-1}^{-1})/(\tau_{n-1}^1\tau_{n+1}^{-1})$,
we have the birational transformations of the form
\begin{align*}
s_n(f_{n-1})&= a_nf_{n-1} \frac{f_n+1/{v_n}}{f_n+u_n},
\quad
s_n(f_{n+1})= \frac{f_{n+1}}{a_{n}} \frac{f_n+u_n}{f_n+1/{v_n}},
\\
r_i(f_n)&= \frac{g_{n+1}(f)}{f_{n+1} g_{n-1}(f)},
\quad
\pi(f_n)=f_{n+1}, 
\quad
\iota(f_n)={1}/{f_n},
\end{align*}
where $g_n(f)=(1+ \sum_{j=1}^{N-1} \prod_{i=1}^{j} u_{n+i}/f_{n+i})
\prod_{k=1}^{N-1}
(u_{n+k})^{-1+k/N}$.
This is in fact equivalent to
the representation due to Kajiwara et al. \cite{kny1, kny2},
which was originally discovered, in connection with a $q$-analogue of the KP hierarchy, without any algebro-geometric setup.
The birational action 
$f_n \mapsto  \overline{f_n} =r_1 \circ \iota(f_n)$ 
of  the translation 
$r_1 \circ \iota:(b_1,b_0;a_n)\mapsto(q b_1,b_0/q;a_n)$ is 
the $q$-Painlev\'e equation of type $A_{N-1}^{(1)}$:
\begin{equation}  \label{eq:q-P(A)}
\overline{f_n}= f_{n+1}\frac{g_{n-1}(f)}{g_{n+1}(f)}
\quad
(n \in {\mathbb Z}/N{\mathbb Z} ).
\end{equation}
The family ${\cal X}=\cup_a X_a$ of 
rational varieties provides 
a sort of phase space for the dynamics of (\ref{eq:q-P(A)}).

\subsection{The case of $D_{N+2}^{(1)}$}
We first need to consider a preliminary situation to 
derive the tropical representation of $W(D_{N+2}^{(1)})$ $(N \geq 3)$;
for instance, let 
${\boldsymbol k}={\boldsymbol \ell}=(2,1,1,\ldots,1,2,1) \in ({\mathbb Z}_{>0})^N$.
Next we shall ignore the simple reflection $s_N^0$ 
corresponding to the $N$-th vertex
and, simultaneously, put $a_N^0=1$ and $\tau_N^{1}=\tau_N^{-1}=1$.
Write the Dynkin diagram of type $D_{N+2}^{(1)}$ as follows:
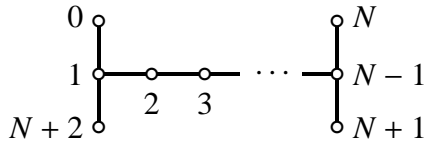
\begin{center}
\begin{picture}(170,50)
\thicklines
\put(40,22){\line(0,1){16}}
\put(40,40){\circle{4}}   \put(28,38){$0$}

\put(40,18){\line(0,-1){16}}
\put(40,0){\circle{4}}   \put(6,-4){$N+2$}

\put(40,20){\circle{4}}     \put(28,16){$1$}
\put(58,20){\line(-1,0){16}}
\put(60,20){\circle{4}}    \put(57,5){$2$}
\put(78,20){\line(-1,0){16}}
\put(80,20){\circle{4}}    \put(77,5){$3$}
\put(82,20){\line(1,0){11}}

\put(99,20){$\ldots$}

\put(128,20){\line(-1,0){11}}
\put(130,20){\circle{4}}    \put(136,16){$N-1$}

\put(130,22){\line(0,1){16}}
\put(130,40){\circle{4}}   \put(136,38){$N$}

\put(130,18){\line(0,-1){16}}
\put(130,0){\circle{4}}   \put(136,-4){$N+1$}

\end{picture}
\end{center}
According to the labels of the diagram above,
we introduce
the multiplicative root variables $a=(a_0, \ldots, a_{N+2})$ and define
the action of simple reflections $s_i$
as 
$s_i(a_j)=a_j {a_i}^{-C_{ij}}$,
where $C_{ij}$ is the Cartan matrix of type $D_{N+2}^{(1)}$.
Set
\[
\begin{array}{l}
u_1=a_1 {a_0}^{-1/2}{a_{N+2}}^{1/2},  \quad v_1=a_1 {a_0}^{1/2}{a_{N+2}}^{-1/2},
\\
u_n=v_n=a_n    \quad (2 \leq n \leq N-2),
\\
u_{N-1}=a_{N-1} {a_{N}}^{-1/2}{a_{N+1}}^{1/2}, \quad
v_{N-1}=a_{N-1} {a_{N}}^{1/2}{a_{N+1}}^{-1/2}.
\end{array}
\]
From Theorem~\ref{thm:tau-rep},
we thus have 
 the following  birational transformations
 $s_n$ $(0 \leq n \leq N+2)$ on the field
 ${\mathbb C}(a^{1/4})(\tau)$ of 
$\tau$-variables:
\begin{align*}
s_n(\tau_n^{1})
&= \frac{{v_n}^{1/2}\xi_{n+1}\eta_{n-1} +{v_n}^{-1/2}\xi_{n-1}\eta_{n+1} }{\tau_n^{-1}},
\\
s_n(\tau_n^{-1})
&= \frac{{u_n}^{-1/2}\xi_{n+1}\eta_{n-1} +{u_n}^{1/2}\xi_{n-1}\eta_{n+1} }{\tau_n^{1}}
\end{align*}
for $1 \leq n \leq N-1$,
and
\[s_0(\tau_1^{\{1,2\}})=\tau_1^{\{2,1\}}, \quad
s_{N+2}(\tau_1^{\{-1,-2\}})=\tau_1^{\{-2,-1\}},
\quad
s_{N}(\tau_{N-1}^{\{1,2\}})
=\tau_{N-1}^{\{2,1\}},
\quad
s_{N+1}(\tau_{N-1}^{\{-1,-2\}})=\tau_{N-1}^{\{-2,-1\}},
\]
where
\[
(\xi_n,\eta_n)= 
\left\{ 
\begin{array}{lll}
(\tau_n^1\tau_n^2,\tau_n^{-1}\tau_n^{-2}) & \text{if} &n=1,N-1,
\\
(\tau_n^1,\tau_n^{-1}) & \text{if}  & 2 \leq n \leq N-2, \\
(1,1) & \text{if}  & n=0,N.
\end{array}
\right.
\]
When $N=3$ ($D_5^{(1)}$ case),  
the above representation in $\tau$-variables coincides with that of \cite[\S1]{tm}.
So, the birational action on $f$-variables
of a translation part of the Weyl group yields 
the $q$-Painlev\'e VI equation \cite{js, s}.
For general $N \geq 3$, 
one may regard it as a (higher order) generalization of the $q$-Painlev\'e VI equation equipped with $W(D_{N+2}^{(1)})$ symmetry.

\begin{remark} \rm
Via the change of variables (\ref{eq:f-tau}), we also have a tropical representation of $W(D_{N+2}^{(1)})$ acting on $f=(f_1, \ldots,f_N)$.
The resulting dynamical system 
looks naively of $N$ dimensional.
However, (i) if $N$ is odd, it possesses a conserved quantity 
$\prod_{n=1}^N f_n$;
(ii) if $N$ is even, it does two 
ones
$\prod_{n=1}^{N/2} f_{2n}$ and $\prod_{n=1}^{N/2} f_{2n-1}$. 
Hence the dimension of the dynamical system is essentially 
$N-1$ (resp. $N-2$) if $N$ is odd (resp. even).
For instance, 
we remember that the $q$-Painlev\'e VI equation ($D_{5}^{(1)}$ case)
is of second order.
\end{remark}

\section{$\tau$-Functions and character polynomials}
\label{sect:char}

In this section, for the case of $A_{N-1}^{(1)}$,
we try to interpret 
from  a soliton-theoretic  point of view
the relationship between 
our $\tau$-functions and character polynomials, i.e.,
the Schur functions or the universal characters.

\subsection{Universal characters}

For a pair of partitions
$\lambda = (\lambda_1,\lambda_2, \ldots,\lambda_r)$
and
$\mu =(\mu_1,\mu_2, \ldots,\mu_{r'}) $, 
the {\it universal character} 
$S_{[\lambda,\mu ]}(x,y)$
is a polynomial in 
$(x,y)=(x_1,x_2,\ldots,y_1,y_2,\ldots)$
defined 
by the determinant formula of
 {\it twisted} Jacobi-Trudi type 
 \cite{k}:
\begin{equation}  \label{eq:def-of-uc}
S_{[\lambda,\mu ]}(x,y)
= \det 
\left(
  \begin{array}{ll}
 p_{\mu_{r'-i+1}  +i - j }(y),  &  1 \leq i \leq r'  \\
 p_{\lambda_{i-r'}-i+j}(x),     &  r'+1 \leq i \leq r+r'   \\
  \end{array}
\right)_{1 \leq i,j \leq r+r'},
\end{equation}
where   
$p_n$ 
is a polynomial defined by 
the generating function
$\sum_{k  \in {\mathbb Z}} p_k(x)z^k =  
\exp 
\left(  \sum_{n=1}^\infty x_n z^n  \right)$
and $p_n(x)=0$ for $n <0$.
The Schur function
$S_\lambda(x)$ 
(see, e.g., \cite{mac}) 
is regarded as
a special case of the universal character: 
\[
S_\lambda(x) 
= \det \bigl( p_{\lambda_i-i+j}(x) \bigr)
= S_{[\lambda, \emptyset ]}(x,y).
\]
If we count the degree of the variables as
$\deg x_n=n$ and
$\deg y_n=-n$,
then 
$S_{[\lambda,\mu ]}(x,y)$ 
is a weighted homogeneous polynomial of degree 
$|\lambda| -|\mu|$,
where $|\lambda|=\lambda_1+\cdots+\lambda_r$.
The universal character $S_{[\lambda,\mu ]}$
describes the irreducible character of a rational representation of 
the general linear group
$GL(n;{\mathbb C})$ corresponding to a pair of partitions $[\lambda,\mu]$,
while the Schur function 
$S_\lambda$
does
that of a polynomial representation
corresponding to a partition $\lambda$;
see \cite{k} for details.

\begin{example}\rm
For $\lambda=(2,1)$ and $\mu=(1)$, we have
\[
S_{[(2,1),(1)]}(x,y)= 
\left|  
 \begin{array}{ccc} 
  p_1(y) & p_0(y)& p_{-1}(y)  \\
  p_1(x) & p_2(x)& p_3(x)  \\
  p_{-1}(x) & p_0(x)& p_1(x)  
 \end{array} 
\right|
=\left( \frac{{x_1}^3}{3}-x_3\right)y_1-{x_1}^2.
\]
\end{example}

\subsection{Notations}
Let us recall some terminology relevant to the case at hand.
A subset ${\bf m} \subset {\mathbb Z}$ is said to be a 
{\it Maya diagram} 
if
$i \in {\bf m}$ (for $i \ll 0$) 
and
$i \notin {\bf m}$ (for $i \gg 0$).
Each Maya diagram 
${\bf m} = \{ \ldots, m_3, m_2, m_1 \}$
corresponds to
a unique partition
$\lambda=(\lambda_1,\lambda_2,\ldots)$ 
such that $m_i-m_{i+1}=\lambda_i-\lambda_{i+1}+1$.
For a sequence of integers
$\nu=(\nu_1,\nu_2,\ldots,\nu_{N})\in{\mathbb Z}^{N}$,
we consider
the Maya diagram
\[
{\bf m}(\nu)=
(N{\mathbb Z}_{<\nu_1}+1)
\cup
(N{\mathbb Z}_{<\nu_2}+2)
\cup
\cdots
\cup
(N{\mathbb Z}_{<\nu_N}+N),
\]
and denote by 
$\lambda(\nu)$
the corresponding partition.
We call a partition of the form $\lambda(\nu)$
an {\it $N$-core partition}.
It is well-known that a partition $\lambda$ is 
$N$-core 
if and only if
$\lambda$ 
has no hook with length of a multiple of $N$.
See \cite{n}.

\begin{example}\rm
If $N=3$ and $\nu=(2,0 ,3)$ then the resulting partition reads $\lambda(\nu)=(4,2,1,1)$.
\end{example}

Also, we prepare the 
following standard notations of $q$-analysis \cite{gr}.

\noindent
{\it $q$-shifted factorials}:
\[
(z;q)_\infty=\prod_{i=0}^{\infty}(1-z q^i),
\quad
(z;p,q)_\infty=\prod_{i,j=0}^{\infty}(1-z p^iq^j),
\]
and
$(z_1,\ldots,z_r;q)_\infty=(z_1;q)_\infty  \cdots (z_r;q)_\infty$,
etc.

\noindent
{\it Elliptic gamma function}:
\[
\Gamma(z;p,q)
= 
\frac{  \left( p q z^{-1};p,q \right)_\infty }{\left(  z;p,q \right)_\infty}.
\]

\subsection{$\tau$-Functions in terms of character polynomials}

In what follows, 
we shall only treat the case of $A_{N-1}^{(1)}$ 
(see \S~4.1). 
Introduce the lattice 
$P(A_{\ell-1})={\mathbb Z}^\ell/{\mathbb Z}({\boldsymbol e}_1+\cdots +{\boldsymbol e}_\ell)$
of rank $\ell-1$,
called the {\it weight lattice} of type $A_{\ell-1}$.
Here
$\{{\boldsymbol e}_i\}$ 
refers to the canonical basis of 
${\mathbb Z}^\ell$. 
Note that $P(A_1)\cong {\mathbb Z}$.
For  an  element 
$\Lambda \in M=W.\{ E_n^\pm\}$
of the Weyl group orbit,
we can associate 
a sequence of integers 
$\nu=(\nu_1,\ldots,\nu_N) \in {\mathbb Z}^N$
and an integer $\kappa \in {\mathbb Z}$
through the formulae
\[
\nu_i-\nu_{i+1}=\langle \Lambda, \check{\alpha}_i \rangle \quad (i \neq N),
\quad
\nu_N-\nu_{1}+1=\langle \Lambda, \check{\alpha}_N \rangle
\quad
{\rm and}  \quad
\kappa=\langle  \Lambda, \check{\beta}_0 \rangle.
\]
Since
$\nu\in {\mathbb Z}^N$ 
is uniquely determined modulo 
${\mathbb Z}({\boldsymbol e}_1+\cdots +{\boldsymbol e}_N)$,
we can regard $(\nu,\kappa)$ as a point of $P(A_{N-1})\times P(A_1)$.
In fact, this is a one-to-one correspondence between 
$M$ and $P(A_{N-1})\times P(A_1)$.
Now let us index the $\tau$ functions as 
\begin{equation}
\tau(\Lambda)=:\sigma_\nu^\kappa.
\end{equation}
For example, we see that 
$\tau(E_n^1)=\sigma_{{\boldsymbol e}_1+\cdots +{\boldsymbol e}_n}^0$
and
$\tau(E_n^{-1})=\sigma_{{\boldsymbol e}_1+\cdots +{\boldsymbol e}_n}^1$.
We have in general
\[
T_n(\sigma_\nu^\kappa)=\sigma_{\nu + {\boldsymbol e}_n}^\kappa
\quad 
(n \in {\mathbb Z}/N{\mathbb Z}),
\quad
\widetilde{T}(\sigma_\nu^\kappa)=\sigma_\nu^{\kappa+1},
\]
where 
$T_n=\pi \circ s_{n+N-2} \circ \cdots \circ s_{n+1} \circ s_n$ 
and $\widetilde{T}=r_1 \circ \iota$ are
the translations of affine Weyl groups
$\widetilde{W}(A_{N-1}^{(1)})$ and $\widetilde{W}(A_{1}^{(1)})$,
respectively.
As previously seen in \S~\ref{sect:tau}, 
the $\tau$-function 
$\sigma_\nu^\kappa$ is expressed as 
a rational function in the indeterminates $\tau_i^1$ and $\tau_i^{-1}$
$(i \in  {\mathbb Z}/N{\mathbb Z})$
by means of the defining polynomial of 
the corresponding divisor.
However, substituting appropriate special values 
for $\tau_i^1$ and $\tau_i^{-1}$,
we have another expression of the $\tau$-function in terms of 
the character polynomials.

\begin{thm}[\rm cf. Kajiwara-Noumi-Yamada \cite{kny2}\bf]
\label{thm:schur}
Define  
\begin{align*}
F(t)
&=\frac{(q^{2N} t^{2N};q^{2N},q^{2N})_\infty}{(q^{2} t^{2};q^{2},q^{2})_\infty}
\Gamma(-q^{N-1}t^{N-1};q^{N-1},q^{N-1}),
\\
H_\nu
&=\prod_{(i,j)\in \lambda(\nu)}(q^{h_{ij}}-q^{-h_{ij}})q^{(i-j)/2},
\end{align*}
where 
$h_{ij}$ is the hook-length, that is, 
$h_{ij}=\lambda_i+\lambda_j'-i-j+1$.
Under the specialization $a_i=q$, $\tau_i^1=F(1/b_0)$ and $\tau_i^{-1}=F(b_1)$
for all $i \in {\mathbb Z}/N {\mathbb Z}$,
the $\tau$-function
$\sigma_{\nu}^\kappa$ is expressed in terms of the Schur function
attached to the $N$-core partition $\lambda(\nu)$
as 
\begin{equation}
\sigma_{\nu}^\kappa
= F(t)H_\nu S_{\lambda(\nu)}(x),
\end{equation}
where $t=q^\kappa/b_0$ and
\[x_n=\frac{t^n+t^{-n}}{n(q^n-q^{-n})}
\quad (n=1,2,\ldots).
\]
\end{thm}

This theorem is based on the principle that
appropriate reductions of infinite-dimensional integrable systems,
such as the KP hierarchy,
yield
Painlev\'e-type equations.
In fact, the Weyl group actions
(\ref{subeq:weyla}) lead us to
the bilinear relation among $\tau$-functions $\sigma^\kappa_\nu \in {\mathbb K}(\tau)$ 
of the form
\begin{equation}  \label{eq:bil-A}
q^{N \nu_i-| \nu | +i-1} 
\left( \prod_{j=1}^{N} \left(\frac{a_{i+j-1}}{q}\right)^{j/N} \right)
\left(t^N- \frac{1}{t^N}\right) 
\sigma^\kappa_\nu
\sigma^\kappa_{\nu + {\boldsymbol e}_i} 
= t \sigma^{\kappa-1}_\nu
\sigma^{\kappa+1}_{\nu + {\boldsymbol e}_i}
-
\frac{1}{t} 
\sigma^{\kappa+1}_\nu
\sigma^{\kappa-1}_{\nu + {\boldsymbol e}_i},
\end{equation}
where
$t=q^\kappa/b_0$.
%and
%$c_i= \prod_{j=1}^{N-1} \left(a_{i+j-1}/q\right)^{-1+j/N}$.
The bilinear relation (\ref{eq:bil-A}) above
is exactly the similarity reduction of 
the $q$-KP hierarchy.
Since the Schur function $S_\lambda$ is a 
homogeneous solution of the $q$-KP hierarchy
and is thus compatible with the similarity reduction,
we can immediately verify Theorem~\ref{thm:schur}.
Moreover, when $N$ is even,  
(\ref{eq:bil-A}) can also be obtained, alternatively, 
from the $q$-UC hierarchy by a similarity reduction \cite{t04, t05a}.
Therefore, the following theorem is a direct consequence of the fact that
the universal character 
$S_{[\lambda, \mu]}$ 
is homogeneous and solves the $q$-UC hierarchy.

\begin{thm}
\label{thm:uc}
In the case  $N=2g+2$ {\rm ($g=1,2,\ldots$)},
define  
\begin{align*}
\widetilde{F}(c,t) 
&=
\frac{(-c q^3 t^2,-c^{-1} q^3 t^2;q^2,q^4)_\infty (q^{4g+4}t^{4g+4};q^{4g+4},q^{4g+4})_\infty}{(q^4 t^{4};q^{4},q^{4})_\infty}
\\
&\times
\Gamma(-c^{1/2} q^{3/2} t, -c^{-1/2} q^{3/2} t; q,q^2)
\Gamma(-q^{2g}t^{2g}; q^{2g},q^{2g}),
\\
\widetilde{H}_\nu
&=
{c_\nu}^{(|\lambda(\nu_{\rm even})|-|\lambda(\nu_{\rm odd})|)/2}
\prod_{(i,j)\in \lambda(\nu_{\rm odd})}(q^{2 h_{ij}}-q^{-2 h_{ij}})
\prod_{(i,j)\in \lambda(\nu_{\rm even})}(q^{-2 h_{ij}}-q^{2 h_{ij}}),
\end{align*}
with
$c_\nu
=c q^{2(| \nu_{\rm even} | - | \nu_{\rm odd} |)}$,
$\nu_{\rm even}=(\nu_2,\nu_4,\ldots, \nu_{2g+2})$,
$\nu_{\rm odd}=(\nu_1,\nu_3,\ldots, \nu_{2g+1})$
and
$c$ being an arbitrary nonzero constant.
Under the specialization $a_{2i+1}=c $,
$a_{2i}=q^2/c$,
$\tau_{2i}^1= \widetilde{F}(c,1/b_0)$, $\tau_{2i}^{-1}=\widetilde{F}(c,b_1)$,
$\tau_{2i+1}^1= \widetilde{F}(c/q^2,1/b_0)$ and $\tau_{2i+1}^{-1}=\widetilde{F}(c/q^2,b_1)$
for all $i \in {\mathbb Z}/(g+1) {\mathbb Z}$,
the $\tau$-function
$\sigma_{\nu}^\kappa$ is expressed in terms of the universal character
attached to a pair of $(g+1)$-core partitions 
as 
\begin{equation}
\sigma_\nu^\kappa
= 
\widetilde{F}(c_\nu,t)  
\widetilde{H}_\nu
S_{[\lambda(\nu_{\rm odd}),
\lambda(\nu_{\rm even})'
]}
(x,y),
\end{equation}
where  $t=q^\kappa/b_0$ and
\[
x_n=\frac{t^{2n}+t^{-2n}-(-{c_\nu})^n(q^n+q^{-n})}{n(q^{2n}-q^{-2n})},
\quad
y_n=\frac{t^{2n}+t^{-2n}-(-{c_\nu})^{-n}(q^n+q^{-n})}{n(q^{-2n}-q^{2n})}
\quad (n=1,2,\ldots).
\]
\end{thm}

Taking into account the origin of our realization of Weyl groups,
Theorems~\ref{thm:schur} and \ref{thm:uc}  
exemplify a quite curious coincidence
between algebraic geometry of particular rational varieties and representation theory. 
Similar results have been established also for the
$D_5^{(1)}$ and $E_6^{(1)}$
cases \cite{t05b, tm} by the use of the $q$-UC hierarchy.
The authors expect there should exist 
another geometric understanding of such a relationship between 
$\tau$-functions and character polynomials,
without employing the theory of integrable systems.

\small
\paragraph{\it Acknowledgement.}
The authors wish to thank Ralph Willox for his helpful comments 
on the first draft of this article.
They have also benefited from discussions
with colleagues in Kobe:
Tetsu Masuda, Masatoshi Noumi,  Yasuhiro Ohta, Masa-Hiko Saito and Yasuhiko Yamada.
This work was supported in part by
a grant-in-aid from the Japan Society for the Promotion of Science (JSPS).

\noindent
Teruhisa Tsuda  \quad
\verb|tudateru@math.kyushu-u.ac.jp|
\\
Faculty of Mathematics, Kyushu University,
Hakozaki, Fukuoka 812-8581, Japan
\\ \\
Tomoyuki Takenawa \quad
\verb|takenawa@kaiyodai.ac.jp|
\\
Faculty of Marine Technology, Tokyo University of Marine Science and Technology,
Echujima, Tokyo 135-8533,
Japan

\end{document}